\documentstyle[twoside]{article}

\catcode`\@=11
\long\def\@makefntext#1{
\protect\noindent \hbox to 3.2pt {\hskip-.9pt
$^{{\eightrm\@thefnmark}}$\hfil}#1\hfill}               

\def\@makefnmark{\hbox to 0pt{$^{\@thefnmark}$\hss}}    

\def\ps@myheadings{\let\@mkboth\@gobbletwo
\def\@oddhead{\hbox{}
\rightmark\hfil\eightrm\thepage}
\def\@oddfoot{}\def\@evenhead{\eightrm\thepage\hfil
\leftmark\hbox{}}\def\@evenfoot{}
\def\sectionmark##1{}\def\subsectionmark##1{}}

\oddsidemargin=\evensidemargin
\newcounter{sectionc}\newcounter{subsectionc}\newcounter{subsubsectionc}
\renewcommand{\section}[1] {\vspace{12pt}\addtocounter{sectionc}{1}
\setcounter{subsectionc}{0}\setcounter{subsubsectionc}{0}\noindent
        {\tenbf\thesectionc. #1}\par\vspace{5pt}}
\renewcommand{\subsection}[1] {\vspace{12pt}\addtocounter{subsectionc}{1}
      \setcounter{subsubsectionc}{0}\noindent
      {\bf\thesectionc.\thesubsectionc.{\kern1pt \bfit #1}}\par\vspace{5pt}}
\renewcommand{\subsubsection}[1]
      {\vspace{12pt}\addtocounter{subsubsectionc}{1}
      \noindent{\tenrm\thesectionc.\thesubsectionc.\thesubsubsectionc.
      {\kern1pt \tenit #1}}\par\vspace{5pt}}
\newcommand{\nonumsection}[1] {\vspace{12pt}\noindent{\tenbf #1}
        \par\vspace{5pt}}

\newcounter{appendixc}
\newcounter{subappendixc}[appendixc]
\newcounter{subsubappendixc}[subappendixc]
\renewcommand{\thesubappendixc}{\Alph{appendixc}.\arabic{subappendixc}}
\renewcommand{\thesubsubappendixc}
        {\Alph{appendixc}.\arabic{subappendixc}.\arabic{subsubappendixc}}

\renewcommand{\appendix}[1] {\vspace{12pt}
        \refstepcounter{appendixc}
        \setcounter{figure}{0}
        \setcounter{table}{0}
        \setcounter{lemma}{0}
        \setcounter{theorem}{0}
        \setcounter{corollary}{0}
        \setcounter{definition}{0}
        \setcounter{equation}{0}
        \renewcommand{\thefigure}{\Alph{appendixc}.\arabic{figure}}
        \renewcommand{\thetable}{\Alph{appendixc}.\arabic{table}}
        \renewcommand{\theappendixc}{\Alph{appendixc}}
        \renewcommand{\thelemma}{\Alph{appendixc}.\arabic{lemma}}
        \renewcommand{\thetheorem}{\Alph{appendixc}.\arabic{theorem}}
        \renewcommand{\thedefinition}{\Alph{appendixc}.\arabic{definition}}
        \renewcommand{\thecorollary}{\Alph{appendixc}.\arabic{corollary}}
        \renewcommand{\theequation}{\Alph{appendixc}.\arabic{equation}}
        \noindent{\tenbf Appendix \theappendixc #1}\par\vspace{5pt}}
\newcommand{\subappendix}[1] {\vspace{12pt}
        \refstepcounter{subappendixc}
        \noindent{\bf Appendix \thesubappendixc. {\kern1pt \bfit #1}}
        \par\vspace{5pt}}
\newcommand{\subsubappendix}[1] {\vspace{12pt}
        \refstepcounter{subsubappendixc}
        \noindent{\rm Appendix \thesubsubappendixc. {\kern1pt \tenit #1}}
        \par\vspace{5pt}}

\newcommand{\smalllineskip}{\baselineskip=10pt}

\def\eightcirc{
\begin{picture}(0,0)
\put(4.4,1.8){\circle{6.5}}
\end{picture}}
\def\eightcopyright{\eightcirc\kern2.7pt\hbox{\eightrm c}}

\def\abstracts#1#2#3{{
        \centering{\begin{minipage}{4.5in}\baselineskip=10pt\footnotesize
        \parindent=0pt #1\par
        \parindent=15pt #2\par
        \parindent=15pt #3
        \end{minipage}}\par}}

\renewenvironment{thebibliography}[1]
        {\frenchspacing
         \ninerm\baselineskip=11pt
         \begin{list}{\arabic{enumi}.}
        {\usecounter{enumi}\setlength{\parsep}{0pt}
         \setlength{\leftmargin 12.7pt}{\rightmargin 0pt} 
         \setlength{\itemsep}{0pt} \settowidth
        {\labelwidth}{#1.}\sloppy}}{\end{list}}

\newcounter{itemlistc}
\newcounter{romanlistc}
\newcounter{alphlistc}
\newcounter{arabiclistc}

\newcommand{\fcaption}[1]{
        \refstepcounter{figure}
        \setbox\@tempboxa = \hbox{\footnotesize Fig.~\thefigure. #1}
        \ifdim \wd\@tempboxa > 5in
           {\begin{center}
        \parbox{5in}{\footnotesize\smalllineskip Fig.~\thefigure. #1}
            \end{center}}
        \else
             {\begin{center}
             {\footnotesize Fig.~\thefigure. #1}
              \end{center}}
        \fi}

\newcommand{\tcaption}[1]{
        \refstepcounter{table}
        \setbox\@tempboxa = \hbox{\footnotesize Table~\thetable. #1}
        \ifdim \wd\@tempboxa > 5in
           {\begin{center}
        \parbox{5in}{\footnotesize\smalllineskip Table~\thetable. #1}
            \end{center}}
        \else
             {\begin{center}
             {\footnotesize Table~\thetable. #1}
              \end{center}}
        \fi}

\def\@citex[#1]#2{\if@filesw\immediate\write\@auxout
        {\string\citation{#2}}\fi
\def\@citea{}\@cite{\@for\@citeb:=#2\do
        {\@citea\def\@citea{,}\@ifundefined
        {b@\@citeb}{{\bf ?}\@warning
        {Citation `\@citeb' on page \thepage \space undefined}}
        {\csname b@\@citeb\endcsname}}}{#1}}

\newif\if@cghi
\def\cite{\@cghitrue\@ifnextchar [{\@tempswatrue
        \@citex}{\@tempswafalse\@citex[]}}
\def\citelow{\@cghifalse\@ifnextchar [{\@tempswatrue
        \@citex}{\@tempswafalse\@citex[]}}
\def\@cite#1#2{{$\null^{#1}$\if@tempswa\typeout
        {IJCGA warning: optional citation argument
        ignored: `#2'} \fi}}

\def\@refcitex[#1]#2{\if@filesw\immediate\write\@auxout
        {\string\citation{#2}}\fi
\def\@citea{}\@refcite{\@for\@citeb:=#2\do
        {\@citea\def\@citea{, }\@ifundefined
        {b@\@citeb}{{\bf ?}\@warning
        {Citation `\@citeb' on page \thepage \space undefined}}
        \hbox{\csname b@\@citeb\endcsname}}}{#1}}

\def\@refcite#1#2{{#1\if@tempswa\typeout
        {IJCGA warning: optional citation argument
        ignored: `#2'} \fi}}

\def\refcite{\@ifnextchar[{\@tempswatrue
        \@refcitex}{\@tempswafalse\@refcitex[]}}

\def\pmb#1{\setbox0=\hbox{#1}
        \kern-.025em\copy0\kern-\wd0
        \kern.05em\copy0\kern-\wd0
        \kern-.025em\raise.0433em\box0}

\def\fnt#1#2{\footnotetext{\kern-.3em
        {$^{\mbox{\scriptsize #1}}$}{#2}}}

\def\fpage#1{\begingroup
\voffset=.3in
\thispagestyle{empty}\begin{table}[b]\centerline{\footnotesize #1}
        \end{table}\endgroup}

\def\runninghead#1#2{\pagestyle{myheadings}
\markboth{{\protect\footnotesize\it{\quad #1}}\hfill}
{\hfill{\protect\footnotesize\it{#2\quad}}}}
\font\tenrm=cmr10
\font\tenit=cmti10
\font\tenbf=cmbx10
\font\bfit=cmbxti10 at 10pt
\font\ninerm=cmr9

\font\eightrm=cmr8

\newtheorem{theorem}{Theorem}

\newtheorem{lemma}{Lemma}

\newtheorem{definition}{Definition}
\newtheorem{corollary}{Corollary}

\newtheorem{example}{Example}

\def\qed{\hbox{${\vcenter{\vbox{                      
   \hrule height 0.4pt\hbox{\vrule width 0.4pt height 6pt
   \kern5pt\vrule width 0.4pt}\hrule height 0.4pt}}}$}}

\begin{document}

\newcommand{\Proof}{{\em{Proof. }}}
\newcommand{\QED}{\ \hfill $\FBox$ \\[1em]}

\newcommand{\M}{\mbox{${\cal M}$}}
\newcommand{\G}{\mbox{${\cal G}$}\ }
\newcommand{\Lsa}{\mbox{$L^{\mbox{\scriptsize sa}} $}}
\newcommand{\Tr}{\mbox{Tr\/}}
\newcommand{\tr}{\mbox{tr\/}}
\newcommand{\Pd}[1]{ \frac{\partial}{\partial x^{#1}} }
\newcommand{\Pdn}[1]{ \frac{\partial}{\partial #1} }
\newcommand{\Pdd}{\mbox{$\partial$ \hspace{-1.2 em} $/$}}
\newcommand{\Op}{{\cal O}\!{\cal P}}
\newcommand{\X}{X}
\newcommand{\Sl}{\mbox{$\prec \!\!$ \nolinebreak}}
\newcommand{\Sr}{\mbox{\nolinebreak $\succ$}}
\newcommand{\bra}{\mbox{$< \!\!$ \nolinebreak}}
\newcommand{\ket}{\mbox{\nolinebreak $>$}}
\newcommand{\kernel}{\mbox{kernel}}
\newcommand{\spc}{\;\;\;\;\;\;\;\;\;\;}
\newcommand{\T}{{\mbox{T }}}
\newcommand{\texp}{{\mbox{Texp }}}
\newcommand{\arctanh}{{\mbox{arctanh}}}
\newcommand{\const}{{\mbox{const }}}

\newcommand{\C}{\mbox{\rm I \hspace{-1.25 em} {\bf C}}}
\newcommand{\R}{\mbox{\rm I \hspace{-.8 em} R}}
\newcommand{\1}{\mbox{\rm 1 \hspace{-1.05 em} 1}}
\newcommand{\Z}{\mbox{\rm \bf Z}}
\newcommand{\sR}{\mbox{\rm \scriptsize I \hspace{-.8 em} R}}
\newcommand{\N}{\mbox{\rm I \hspace{-.8 em} N}}
\newcommand{\sN}{\mbox{\rm \scriptsize I \hspace{-.8 em} N}}
\newcommand{\sZ}{\mbox{\rm \scriptsize \bf Z}}
\newcommand{\loc}{\mbox{\rm{\scriptsize{loc}}}}
\newcommand{\nab}{\breve{\nabla}}
\newcommand{\rhe}{\rho_{\mbox{\rm \scriptsize em}}}
\newcommand{\rhm}{\rho_{\mbox{\rm \scriptsize m}}}
\newcommand{\sj}{\breve{\jmath}}
\newcommand{\srhe}{\breve{\rho}_{\mbox{\rm \scriptsize em}}}
\newcommand{\srhm}{\breve{\rho}_{\mbox{\rm \scriptsize m}}}

\newcommand{\Equ}[1]{\begin{equation} \label{eq:#1}}
\newcommand{\EndEqu}{\end{equation}}
\newcommand{\Ref}[1]{(\ref{eq:#1})}

\newcommand{\FBox}{\rule{2mm}{2.25mm}}
\newcommand{\OBox}{\raisebox{.6ex}{\fbox{}}\,}
\newcommand{\A}{{\mbox{${\cal A}$}}}
\newcommand{\B}{{\mbox{${\cal B}$}}}

\runninghead{F. Finster}
{Derivation of Local Gauge Freedom $\ldots$}

\thispagestyle{empty}\setcounter{page}{315}
\vspace*{0.88truein}
\fpage{315}

\centerline{\bf DERIVATION OF LOCAL GAUGE FREEDOM}
\centerline{\bf FROM A MEASUREMENT PRINCIPLE}

\vspace*{0.035truein}

\vspace*{0.37truein}
\centerline{\footnotesize Felix Finster}

\centerline{\footnotesize \it
Department of Mathematics, Harvard University}
\baselineskip=10pt
\centerline{\footnotesize \it
E-mail: felix@math.harvard.edu}

\baselineskip 5mm

\vspace*{0.21truein}

\abstracts{ We define operator manifolds as manifolds on which a spectral
measure on a Hilbert space is given as additional structure.  The spectral
measure mathematically describes space as a quantum mechanical observable.
We show that the vectors of the Hilbert space can be represented as
functions on the manifold. The arbitrariness of this representation is
interpreted as local gauge freedom. In this way, the physical gauge principle
is linked with quantum mechanical measurements of the position variable.
We derive the restriction for the local gauge group to be $U(m)$, where $m$
is the number of components of the wave functions.}{}{}

\bigskip

$$$$

\section{\bf Operator Manifolds}
In nonrelativistic quantum mechanics, the physical observables ${\cal{O}}$
are described by self-adjoint operators on a Hilbert space $H$.
The wave function of a particle is a vector $\Psi \in H$; measurements
correspond to the calculation of expectation values $\bra \Psi |
{\cal{O}} | \Psi \ket$.
In this paper, we will concentrate on the observables for space.
These observables play a special role because they
determine the geometry of the physical system.
Usually they are given by mutually commuting operators $(X^i)_{i=1,\ldots,3}$
with a continuous spectrum (i.e.\ the multiplication operators with the
coordinate functions in position space).
Since the $X^i$ depend on the choice of the coordinate system
in $\R^3$, it is more convenient to consider their spectral
measure $(E_x)_{x \in \sR^3}$ (for basic definitions
see e.g.\ \cite{reed}). The operators $X^i$ can be reconstructed from the
spectral measure by integrating over the coordinate functions,
\Equ{1}
X^i \;=\; \int_{\sR^3} x^i \: dE_x \;\;\; .
\EndEqu
We want to study this functional analytic point of view in the more general
setting that space is a manifold.
\begin{definition}
Let $M$ be a manifold of dimension $n$, $(\mu, {\cal{M}})$ a positive,
$\sigma$-finite measure on $M$ with $\sigma$-algebra ${\cal{M}}$.
Furthermore, let $H$ be a separable Hilbert space,
\[ E \;:\; {\cal{M}} \;\rightarrow\; P(H) \]
a spectral measure ($P(H)$ denotes the projection operators of $H$), which is
absolutely continuous with respect to $\mu$, $E \ll \mu$
(for basic measure theory see e.g. \cite{rudin}).
$(M, \mu, H, E)$ is called {\bf{operator manifold}}.
\end{definition}
For simplicity, the reader may think of $d\mu$ as the measure
$\sqrt{g} \: d^nx$ on a Riemannian manifold and of ${\cal{M}}$ as the Lebesgue
measurable sets. The requirement $E \ll \mu$ is mainly a technical
simplification.

\begin{example}
\begin{enumerate}
\item \em For a scalar particle, we choose $M=\R^3$, $H=L^2(M)$, and
$(\mu, {\cal{M}})$ the Lebesgue measure. We define the projectors
\[ P_V \;:\; H \:\rightarrow\: H \;:\; f \:\rightarrow\: \chi_V \: f
	\;\;\;,\spc V \in {\cal{M}} \]
as the multiplication operators with the characteristic function.
For the spectral measure $E(V)=P_V$, the integrals \Ref{1} give the usual position operators of quantum mechanics.
\item For a particle with spin $\frac{1}{2}$, we choose $M=\R^3$,
$(\mu, {\cal{M}})$ the Lebesgue measure, and $H=L^2(M) \oplus L^2(M)$
the two-component spinors.
We set $P_V(f^\alpha) = \chi_V f^\alpha$ and again define the spectral measure
by $E(V)=P_V$.
\item In order to describe a scalar particle whose motion is (by some external
forces or constraints) restricted to a submanifold $M \subset \R^3$, we
take $H=L^2(M)$, $(\mu, {\cal{M}})$ the measure $\sqrt{g} \: d^nx$ on $M$
(for the induced Riemannian metric), and
 $E(V)=P_V$, $P_V(f)=\chi_V f$. For a chart $(x^i, U)$, the
corresponding position operators are given by
\[ X^i \;=\; \int_U x^i \; dE_x \;\;\; . \] \em
\end{enumerate}
\end{example}
In these examples, the vectors of $H$ are functions on the manifold with one or
several components. We want to study the question if $H$ can also in the general
case be represented as a function space on $M$.
A possible method for this analysis is the functional calculus and
constructions similar to the proof of the spectral theorem in its multiplicative
form (see e.g. \cite{reed}).
We proceed in a different way using the notions of ``spin scalar product''
and ``local orthonormal basis,'' which is considered to be more transparent.

\begin{definition}
For $u,v \in H$, we define the complex, bounded measure $\mu_{uv}$ by
\[ \mu_{uv}(V) \;=\; \bra E_V u, v \ket \;\;\; . \]
Since $\mu_{uv} \ll \mu$, it has a
unique Radon-Nikodym representation $d\mu_{uv}=h_{uv} \: d\mu$
with $h_{uv} \in L^1(M, \mu)$. The mapping
\[ \Sl .,. \Sr \;:\; H \times H \:\rightarrow\: L^1(M, \mu) \;:\;
	(u,v) \:\rightarrow\: h_{uv} \]
is called {\bf{spin scalar product}}.
\end{definition}
The spin scalar product is linear in the first and anti-linear in
the second argument. Furthermore it is positive, $\Sl u,u \Sr \geq 0$.

\begin{definition}
\label{def_ONB}
A family $(u_l, C_l)_{l \in \sN}$ with $u_l \in H$, $C_l \in {\cal{M}}$ is
called {\bf{local orthonormal basis}} (local ONB) if
\begin{description}
\item[(i)] $\Sl u_l,\: u_m \Sr \;=\; \delta_{lm} \: \chi_{C_l}$
\item[(ii)] The subspace $\bra \left\{ E_V u_l \::\: V \in {\cal{M}},\: l \in
	\N \right\} \ket$ is dense in $H$.
\end{description}
We define a measurable partition $(D_m)_{m \in \sN \cup \{0, \infty\}}$ of $M$ by
\Equ{Y}
D_m \;=\; \left\{ x \in M \:|\: \# \{ l \:|\: x \in C_l \}=m \right\} \;\;\; .
\EndEqu
We say that on $D_m$ the {\bf{spin dimension}} is $m$.
\end{definition}

\begin{lemma}
    \label{ex_loc_ONB}
    There exists a local orthonormal basis.
\end{lemma}
{\Proof}
Let $(v_{l})_{l \in \sN}$ be an orthonormal basis of $H$. We construct
a local ONB $(u_l, C_l)$ in several steps:
\begin{enumerate}
    \item Using the notation
\Equ{S}
H_u \;=\; \overline{ \bra \{ E_V u : V \in \M \} \ket } \;\;\;,\spc
	u \in H \;\;\; ,
\EndEqu
we define the closed subspaces $K_l$ by $K_l = \bra H_{v1}, \ldots, H_{v_{l}} \ket$
and construct the series $(w_l)_{l \in \sN}$ by
$w_1=v_1$, $w_l = (1- \Pr_{K_{l-1}})\: v_l$
($\Pr_K$ denotes the projector on the closed subspace $K$).
The $K_l$ are a filtration of $H$, i.e. $K_l \subset K_{l+1}$ and $\overline{
\cup_l K_l} = H$. Furthermore, $K_{l-1} \perp w_l \in K_l$ and
$K_l = \bra H_{w_1}, \ldots, H_{w_l} \ket$.
For $l<m$, $\bra E_V w_l, E_W w_m \ket = \bra E_{V \cap W} w_l,w_m \ket = 0$
and thus $H_{w_l} \perp H_{w_m}$.
We conclude that $H = \bigoplus_l H_{w_l}$ and $\Sl w_l, w_m \Sr = 0$
for $l \neq m$.
    \item
        Set $C_{l} = \{ x \:|\: \Sl w_{l},w_{l} \Sr(x) \neq 0 \}$.
	We may assume that $\mu(C_{l}) < \infty$ for all $l$, because we can
        otherwise take a partition $(U_{i})_{i \in \sN} \in \M$ of $M$
        with $\mu(U_{i}) < \infty$\ and define the vectors $w_{li} = E_{U_{i}} w_{l}$. Then the sets $C_{li} = \{ x \:|\: \Sl w_{li},w_{li} \Sr \neq 0 \}
        = C_{l} \cap U_{i}$\  have finite measure.
        Since in addition $H_{w_l} = \bigoplus_i H_{w_{li}}$,
	we can replace the series $(w_l)$ by $(w_{li})$.
     \item  The functions
\[      f_l \;=\; \frac{\chi_{C_l}}{\sqrt{\Sl w_l,w_l \Sr}}            \]
        are in $L^{2}(M,\mu_{w_l})$, because
\[        \int_{M} |f_l|^2 \:d\mu_{w_l}
            \;=\; \int_M |f_l|^2 \:\Sl w_l,w_l \Sr \:d\mu
	\;=\; \mu(C_l) \;<\; \infty \;\;\; . \]
        Acoording to the spectral theorem for unbounded self-adjoint
        operators, we can thus introduce the vectors $u_l$ by
\[        u_l = \left( \int_{M} f_l(x) \: dE_x \right) \: w_l \;\;\; .   \]
        They satisfy the equation
\[      \int_{V} \Sl u_{l},u_{m} \Sr \: d\mu
                \;=\; \int_{V} d\bra E_{x} u_{l},u_{m} \ket
		\;=\; \int_{V} f_l \: \overline{f_m} \: \Sl w_{l},w_{m} \Sr \: d\mu
		\;\;\; . \]
	In the case $l \neq m$, we obtain $\Sl u_{l},u_{m} \Sr = 0$.
	For $l=m$, we get
\[       \int_V \Sl u_{l},u_{l} \Sr \: d\mu
             = \int_{V} \frac{\chi_{C_{l}}}{\Sl w_{l},w_{l} \Sr}
         \: \Sl w_{l},w_{l} \Sr \: d\mu
             = \int_{V} \chi_{C_{l}} \: d\mu    \;\;\; ,                  \]
and hence $\Sl u_{l},u_{l} \Sr = \chi_{C_{l}}$.
We conclude that $\Sl u_l,u_m \Sr = \delta_{lm} \chi_{C_l}$.
Since
\[         H_{u_l} = \left\{ (\int_M f(x) \: dE_x) \:u_l \;,\; \ f \in
             L^2(M,\mu_{u_l}) \right\} = H_{w_l} \;\;\; ,                 \]
we have $H = \bigoplus_l H_{u_l}$. Thus condition (ii) in Definition~\ref{def_ONB}
is satisfied. \QED
\end{enumerate}
Let in the following $(u_l, C_l)$ be a given local ONB.
\begin{lemma}
    \label{lemma_unit}
The functions $\Sl v,u_l \Sr$ are in $L^2(C_l, \mu)$. The mapping
\Equ{Z}
	U \::\: H \rightarrow \bigoplus_l L^2(C_l,\mu) \::\:
                v \rightarrow \Sl v,u_l \Sr
\EndEqu
      is unitary and $U E_V U^{-1} = P_V$, where $P_V : f_l
\rightarrow f_l \:\chi_V$ is the multiplication operator with the characteristic function.
\end{lemma}
{\Proof}
Using the notation \Ref{S}, Definition~\ref{def_ONB} implies that
\Equ{dir_sum}      H = \bigoplus_l H_{u_l}    \;\;\; .      \EndEqu
We proceed in several steps:
\begin{enumerate}
\item We define the operators
\[      A_l \;:\; <\{E_V u_l : V \in \M \}> \: \rightarrow L^2(C_l,\mu)
	\spc{\mbox{by}}\spc
        A_l (E_V u_l) = \chi_{(V \cap C_l)}. \]
They are isometric, because
\begin{eqnarray*}
\bra E_V u_l, E_W u_l \ket
&=& \int_M \chi_{(V \cap W)} \: \Sl u_l,u_l \Sr \:d\mu \\
&=& \mu(V \cap W \cap C_l) \;=\; \bra \chi_{V \cap C_l}, \chi_{W \cap C_l}
	\ket_{L^2(C_l)}  \;\;\; .
\end{eqnarray*}
       Since ${\cal{D}}(A_l)$ is dense in $H_{u_l}$ and
       ${\cal{R}}(A_l)$ is dense in $L^2(C_l)$, the $A_l$ can
       be uniquely extended to unitary operators
       $A_l : H_{u_l} \rightarrow L^2(C_l)$.
Using the decomposition \Ref{dir_sum}, we define a unitary operator $A$ by
\Equ{u}
	A \;=\; \bigoplus_l A_l \;\;\; .
\EndEqu
\item For $V \in \M$ and $l \in \N$, the vector $u:=E_v \:u_l$ satisfies the
equation
\begin{eqnarray*}
A^\ast P_W A \:u &=& A^\ast P_W (\chi_{(V \cap C_l)} \: \delta_{lm})_{m \in \sN}
        \;=\; A^\ast (\chi_{(W \cap V \cap C_l)} \: \delta_{lm})_{m \in \sN} \\
&=& E_{W \cap V} \: u_l = E_W \: u \;\;\; .
\end{eqnarray*}
Thus $E_W \:u = A^\ast P_W A \:u$ on $\{ E_V u_l : V \in \M \}$.
By continuity, this equation also holds on $H_{u_l}$.
Again by continuity and \Ref{dir_sum}, we obtain $E_W = A^\ast P_W A$.
\item We want to prove that $A=U$. Notice that it is not sufficient
        to show that $A=U$ on a dense subset of $H$, because the continuity
	of $U$ is not obvious.
        Therefore let $v \in H$ be an arbitrary vector and set
            $(f_l)_{l \in \sN} = Av$.
	Since $f_l$ vanishes outside $C_l$, we have
\begin{eqnarray*}
\int_W (Uv)_l \: d\mu &=& \int_W \Sl v,u_l \Sr \: d\mu
	\;=\; \bra E_W v,u_l \ket \\
&=& \bra E_W A^\ast (f_m), u_l \ket
               \;=\; \bra A^\ast A E_W A^\ast (f_m), u_l \ket \\
&=& \bra A^\ast P_W (f_m), u_l \ket \;=\; \bra P_W (f_m), A u_l \ket
                _{\bigoplus L^2(C_l)} \\
&=& \int_{W \cap C_l} f_l \: d\mu \;=\; \int_W f_l \: d\mu \;\;\; .
\end{eqnarray*}
        It follows that $(Uv)_l = f_l = (Av)_l$ and hence $U=A$. \QED
\end{enumerate}
The unitary operator \Ref{Z} gives the desired representation of the
vectors of $H$ as functions on $M$.
The representation is not unique; it depends on the choice of the
local ONB.

The formalism of local ONBs has some analogy with the representation
$u = \sum_l \bra u, u_l \ket \: u_l$ of a vector in an orthonormal basis
$(u_l)$. As the main difference, the scalar product can be ``localized'' on operator
manifolds with the spectral measure, leading to $L^2$-component functions
$\Sl u, u_l \Sr$ instead of complex coefficients $\bra u, u_l \ket$.
The following corollary extends the formal analogy between
ONBs and local ONBs to Parseval's equations.

\begin{corollary}[local completeness relation]
    \label{loc_compl_rel}
    For $u,v \in H$,
\begin{eqnarray}
\label{eq:U}
u &=& \sum_{l=1}^\infty \;\left( \: \int_{C_l}
            \Sl u,u_l \Sr_{|x} \: dE_x \right) \: u_l \\
\label{eq:V}
\Sl u,v \Sr &=& \sum_{l=1}^\infty \Sl u,u_l \Sr \:
            \Sl u_l,v \Sr   \spc    {\mbox{a.e.}} \;\;\; .
\end{eqnarray}
\end{corollary}
{\Proof}
According to Lemma~\ref{lemma_unit}, the function $\Sl u,u_l \Sr$ is in
        $L^2(C_l,\mu)=L^2(C_l,\mu_{u_l u_l})$. We can thus apply the spectral
        theorem for unbounded, self-adjoint operators and define the vectors
\[        w_l = ( \int_{C_l} \Sl u,u_l \Sr_{|x} \: dE_x ) \: u_l \;\;\; .   \]
We have $\Sl w_l, u_m \Sr = \Sl u, u_l \Sr \: \delta_{lm}$ and thus, with
the notation \Ref{S}, $U w_l = U \Pr_{H_{u_l}} u$. The injectivity of $U$ and
\Ref{dir_sum} yield equation \Ref{U}.

        Since $\Sl u,u_l \Sr , \Sl v,u_l \Sr \in \bigoplus_l L^2(C_l)$,
	Lebesgue's monotone convergence theorem yields that
\[        \infty \;>\; \sum_l \int_M | \Sl u,u_l \Sr \:
                \Sl u_l,v \Sr | \: d\mu
          \;=\; \int_M \sum_l | \Sl u,u_l \Sr \: \Sl u_l,v \Sr | \: d\mu
	\;\;\; , \]
	so that the function $f:=\sum_l |\Sl u,u_l \Sr \Sl u_l,v \Sr|$
	is in $L^1(M, \mu)$.
	According to Lemma~\ref{lemma_unit},
\begin{eqnarray}
	\int_W \Sl u,v \Sr \: d\mu &=& \bra E_W u,v \ket \;=\;
            \bra U E_W U^\ast \: U u, Uv \ket_{\bigoplus_l L^2(C_l)} \nonumber \\
\label{eq:t}
	&=& \bra  P_W Uu,Uv \ket \;=\; \sum_l \int_W \Sl u,u_l \Sr \:
	    \Sl u_l,v \Sr \:d\mu \;\;\; .
\end{eqnarray}
        Since the function $f$ dominates the integrand, we can reverse
        the order of summation and integration in \Ref{t} and obtain \Ref{V}.
\QED
We come to the question of how the representation of $H$ as a space of functions
on $M$ depends on the choice of the local ONB. We start with a technical lemma.

\begin{lemma}
    \label{m_eq_n}
    Let $W$ be a measurable set with $\mu(W) \neq 0$,
\[      U \::\: (L^2(W,\mu))^m  \rightarrow (L^2(W,\mu))^n
	\;\;\;\;{\mbox{for}}\;\;\;\ m, n \in \N \cup \{0, \infty\}            \]
    a unitary operator satisfying $U P_V U^{-1} = P_V$ for all
    $V\in \M$.
    Then $m=n$.
\end{lemma}
{\Proof}
    We can assume that $\mu(W) < \infty$, because otherwise take
    a measurable set $V \subset W$ with $0 \neq \mu(V) < \infty$ and
    consider the restriction of $U$ on $(L^2(V,\mu))^m$.
    Since the roles of $m, n$ can be interchanged, it suffices
    to prove $m \leq n$.

    Assume that $m>n$. The vectors
\[      u^\alpha \;:=\; (\delta^\alpha_i)_{i=1,\ldots,m} \in (L^2(W,\mu))^m
	\;\;\;,\spc \alpha=1,\ldots,n+1 \]
satisfy the equations
\[      \bra  P_V \:u^\alpha, u^\beta  \ket = \delta^{\alpha \beta} \mu(V)
        \;\;\;\; {\mbox{ for $V \subset W$}} \;\;\; .                         \]
    By hypothesis on $U$, the functions
    $v^\alpha := U(u^\alpha) = (g^\alpha_i)_{i=1,\ldots,n}$ also satisfy
\[      \bra  P_V \: v^\alpha, v^\beta  \ket = \delta^{\alpha \beta} \mu(V)
        \;\;\;\; {\mbox{ for $V \subset W$}} \;\;\; .                  \]
We thus have for any measurable set $V \subset W$ with $\mu(V) \neq 0$,
\[      \delta^{\alpha \beta}
        = \frac{1}{\mu(V)} \:\bra P_V v^\alpha, v^\beta  \ket
        = \frac{1}{\mu(V)} \int_V \sum_{j=1}^n g^\alpha_j \:
            \overline{g^\beta_j} d\mu  \;\;\; , \]
and thus
\[      \sum_{j=1}^n g^\alpha_j \:
            \overline{g^\beta_j} = \delta^{\alpha \beta}
        \;\; {\mbox{a.e.}} \;\;\;,\spc \alpha, \beta = 1,\ldots, m \;\;\; . \]
This is a contradiction because there are at most $n$ linearly independent
    vectors in $\C^n$.
\QED

\begin{definition}
\label{def_iso}
    Two operator manifolds
    $(M,\mu,H,E), (\tilde{M},\tilde{\mu},\tilde{H},\tilde{E})$
    are called {\bf isomorphic} if
    there exists a homeomorphism $\phi : M \rightarrow \tilde{M}$\
    and a unitary operator $U : H \rightarrow \tilde{H}$ such that
    \begin{description}
      \item[(i)] $\phi$\ preserves the measure:
\[          \mu(V) = \tilde{\mu}(\phi(V)) \;\;{\mbox{for all $V \in \M$}} \]
      \item[(ii)] the spectral measure is invariant:
\[          U E_V U^{-1} = \tilde{E}_{\phi (V)}
                    \;\;{\mbox{for all $V \in M$}} \]
    \end{description}
    In the special case $M=\tilde{M}$ and $\phi = 1$, the operator $U$
    satisfies
\[          U E_V U^{-1} = \tilde{E}_V  \;\;{\mbox{for all $V \in \M$}}   \]
    and is called {\bf{isomorphism}}.
\end{definition}

\begin{theorem}
    \label{class_op_man}
    Let $(M,\mu,H,E)$ be an operator manifold.
    \begin{enumerate}
      \item There exists a partition of M of measurable sets
        $(D_m)_{m \in \sN \cup \{0, \infty\} }$ such that
        $(M,\mu,H,E)$ is isomorphic to the operator manifold
\[       (M, \mu, \bigoplus_{m \in \sN \cup \{0, \infty\}  }
            (L^2(D_m,\mu))^m, P) \;\;\; ,                                     \]
        where $P$ is the canonical spectral measure
        $P : V \rightarrow P_V$ with $P_V : f_i \rightarrow f_i \, \chi_V$.
        The sets $D_m$ are unique (modulo sets of measure zero).
	They coincide with the partition $(D_m)$ in Definition~\ref{def_ONB},
	so that the definition of the spin dimension is independent of the
        choice of the local ONB.
      \item We denote the spin dimension at $x$ by $m_x$ (so $m_x=m$ for $x \in D_m$).
        For two isomorphisms $V,\tilde{V} : H \rightarrow
        \bigoplus_{m \in J  } (L^2(D_m,\mu))^m$, there are measurable functions
	$W^\alpha_\beta$ with
\Equ{loc_unit1}
          \sum_{\alpha=1}^{m_x} W^\alpha_\beta(x) \:
            \overline{W^\alpha_\gamma(x)}
            = \delta_{\beta \gamma} \;\;\;,\spc
	\sum_{\alpha=1}^{m_x} W^\beta_\alpha(x)\:
            \overline{W^\gamma_\alpha(x)}
            = \delta^{\beta \gamma} \spc {\mbox{for a.a. x}},
\EndEqu
such that
\Equ{repr}
         ((V \tilde{V}^{-1}) \: f)^\alpha_{|x} \;=\; \sum_{\beta=1}^{m_x}
                W^\alpha_\beta(x) \: f^\beta(x) \;\;\; .
\EndEqu
        Conversely, if $V$ is an isomphism and $W^\alpha_\beta(x)$ a family of
        measurable functions with the property \Ref{loc_unit1}, then there is
        an isomorphism $\tilde{V}$ satisfying \Ref{repr}.
      \item Every isomorphism $V : H \rightarrow \bigoplus_m (L^2(D_m, \mu))^m$
	can, for a suitable local ONB $(u_l, C_l)$, be realized as the mapping
	$\Sl .,u_l \Sr$.
    \end{enumerate}
\end{theorem}
{\Proof}
    \begin{enumerate}
\item We take a local ONB $(u_l, C_l)$ and define the sets $D_m$ by \Ref{Y}.
	Lemma~\ref{lemma_unit} gives an isomorphism $U$ of $(M, \mu, H, E)$ and
	$(M, \mu, \oplus_l L^2(C_l, \mu), P)$. By cutting and recomposing the
	component functions, we construct a unitary transformation
\Equ{T}
W \::\: \bigoplus_l L^2(C_l, \mu) \;\rightarrow\; \bigoplus_{m \in \sN \cup
	\{0, \infty\} } (L^2(D_m, \mu))^m
\EndEqu
with $W P_V W^{-1} = P_V$. The operator $W U \::\:
	H \rightarrow \bigoplus_m (L^2(D_m, \mu))^m$ is the desired isomorphism.

      In order to show the uniqueness of the sets $D_m$, we consider two
	  isomorphisms
\begin{eqnarray*}
V : H &\rightarrow&
         \bigoplus_{m \in \sN \cup \{0, \infty\} } (L^2(D_m,\mu))^m          \\
\tilde{V} : H &\rightarrow&
         \bigoplus_{m \in \sN \cup \{0, \infty\} } (L^2(\tilde{D}_m,\mu))^m
\end{eqnarray*}
      constructed from different local ONBs. Then the mapping
\[ \tilde{V} V^{-1} :
        \bigoplus_{m} (L^2(D_m,\mu))^m \;\rightarrow\;
        \bigoplus_{m} (L^2(\tilde{D}_m,\mu))^m\]
is an isomorphism. Assume that $\mu(D_m \cap \tilde{D}_m) \neq \mu(D_m)$ or
$\mu(D_m \cap \tilde{D}_m) \neq \mu(\tilde{D}_m)$ for certain $m$.
Since the roles of $D_m$ and $\tilde{D}_m$ can be exchanged, we can assume
that there is a set $W$, $\mu(W) \neq 0$ with $W \subset D_m, W \subset \tilde{D}_n$
and $n \neq m$. Then the restriction of $\tilde{V} V^{-1}$ to
      $P_W (\bigoplus_{m \in J} (L^2(D_m,\mu))^m)$
      is a unitary mapping
\[      A : (L^2(W,\mu))^m \rightarrow (L^2(W,\mu))^n                     \]
      with $A P_V A^{-1} = P_V$ for all $V \subset W$.
      This is a contradiction to Lemma~\ref{m_eq_n}.
      We conclude that $D_m, \tilde{D}_m$ coincide up to sets of measure zero.
     \item We set $U = V \tilde{V}^{-1}$. For $u\!=\!(f^\alpha),v\!=\!(g^\alpha)
      \in \bigoplus_m (L^2(D_m,\mu))^m$ and $C \in \M$,
\[      \int_C \sum_{\alpha=1}^{m_x} f^\alpha \:
            \overline{g^\alpha} d\mu
        \;=\; \bra  P_C \:u,v  \ket \;=\; \bra  P_C \:Uu,Uv  \ket
        \;=\; \int_C \sum_{\alpha =1}^{m_x} (Uu)^\alpha \:
            \overline{(Uv)^\alpha} d\mu           \]
      (integration and summation can be exchanged according to
      Lebesgue's dominated convergence theorem). Hence
\Equ{l_unit}
        \sum_{\alpha=1}^{m_x} f^\alpha(x) \:\overline{g^\alpha(x)}
            \;=\; \sum_{\alpha=1}^{m_x} (Uu)^\alpha(x)\:
            \overline{(Uv)^\alpha(x)}
            \;\;\; {\mbox{a.e.}} \;\;\; ,
\EndEqu
      and we conclude that
\[      (Uu)^\alpha(x) = \sum_{\beta=1}^{m_x}
            U^\alpha_\beta(x) \: f^\beta(x) \;\;\; {\mbox{a.e.}} \]
      for suitable coefficients $U^\alpha_\beta(x)$, which are
      measurable functions in $x$. The identities \Ref{loc_unit1} follow from
	  \Ref{l_unit} and the unitarity of $U$.
      Conversely, if some measurable functions
      $(U^\alpha_\beta(x))_{\alpha,\beta=1,\ldots,m_x}$
      satisfy~\Ref{loc_unit1}, we define the isomorphism $\tilde{V}$ by
\[      \tilde{V} : u \rightarrow ( \:
            \sum_{\gamma=1}^{m_x} \overline{U^\gamma_\beta(x)}
                \: (Vu)^\gamma \: )^\beta \;\;\; .     \]
     \item Let $V : H \rightarrow \bigoplus_m (L^2(D_m, \mu))^m$ be an isomorphism.
	We choose a partition $(C_l)_{l \in \sN}$ of $M$
      subordinate to $D_m$ with $\mu(C_l) < \infty$ and define the mapping
      $m_l : \N \cup \{0, \infty\} \rightarrow \N \cup \{0, \infty\}$
      by the requirement that $m_l = m$ if $C_l \subset D_m$. The vectors
\[      (u_{l \alpha})_{l \in \sN, \alpha=1,\ldots,m_l}
          \;:=\; U^{-1} \: (\chi_{C_l} (\delta^{\alpha \beta}))^\beta      \]
      satisfy for every $W \in \M$ the relations
\begin{eqnarray*}
      \int_W \Sl u_{k \alpha},u_{l \beta} \Sr \: d\mu
         &=& \bra  E_W u_{k \alpha},u_{l \beta}  \ket
         \;=\; \bra  P_W U u_{k \alpha},U u_{l \beta}  \ket \\
&=& \int_W \chi_{C_k} \chi_{C_l} \delta_{\alpha \beta} d\mu
         \;=\; \delta_{\alpha \beta} \delta_{kl} \: \mu(C_k \cap W)  \;\;\; ,
\end{eqnarray*}
      and thus $\Sl u_{k \alpha},u_{l \beta} \Sr = \delta_{\alpha \beta}
            \delta_{kl} \: \chi_{C_k}$.
      Since in addition
\[      H_{u_{k \alpha}} = U^{-1} ( \iota_\alpha L^2(C_k))        \]
      (with $\iota_\alpha : L^2(C_k) \hookrightarrow \sum_m L^2(D_m)^m$ the
	natural injection), it follows that
      $H = \bigoplus_{l,\alpha} H_{u_{l \alpha}}$.
      We conclude that $(u_{l \alpha}, C_l)$ is a local ONB.
By construction of $(u_{l \alpha}, C_l)$, we have $V = W \circ
	(\Sl .,u_l \Sr)_{l \in \sN}$, where $W$ is the canonical isomorphism
	\Ref{T}.
\QED
    \end{enumerate}

\section{\bf Interpretation, Local Gauge Transformations}
In order to clarify the concept of measurability of space, we defined
operator manifolds in a general and abstract mathematical setting.
For physical applications, one needs to introduce additional
objects like the Hamiltonian and the physical wave functions.
Four our discussion, however, this is not necessary; we prefer to keep the
physical interpretation on a general level.

We begin with the case $D_1=M$ of spin dimension one. According to Theorem
\ref{class_op_man}, $H$ is isomorphic to $L^2(M)$ and can be interpreted as the
configuration space of a scalar particle. In contrast to Example 1.2.1, where
$H=L^2(M)$ coincided with the function space, $H$ now is only an abstract
Hilbert space. The isomorphisms of Theorem \ref{class_op_man} give
representations of the vectors of $H$ as ``wave functions.''
The arbitrariness \Ref{loc_unit1},\Ref{repr} of the representation describes
local phase transformations
\Equ{k1}
\Psi(x) \;\longrightarrow\; e^{i \Lambda(x)} \: \Psi(x)
\EndEqu
of the wave functions. This result can be understood qualitatively from the
fact that the wave function itself is not observable, only its absolute square
$|\Psi(x)|^2$ has a physical interpretation as probability density of
the particle.
The transformation \Ref{k1} occurs in quantum mechanics as the local
$U(1)$-gauge transformation of the magnetic field (under which the vector
potential behaves like $\vec{A} \rightarrow \vec{A} + (\vec{\nabla} \Lambda)$).

In the case $D_2=M$ of spin dimension 2, $H$ is isomorphic to $L^2(M) \oplus L^2(M)$
and can be interpreted as two-component Pauli spinors. According to
\Ref{loc_unit1},\Ref{repr}, the arbitrariness of the representation as wave
functions now describes local $U(2)$-transformations.
These transformations really occur in physics; they correspond to the
local $U(2)$-symmetry in quantum mechanics \cite{Fr}.

The general case allows for the description of $m$-component wave functions,
which is needed for particles with higher spin and for particle multipletts.
We again interpret the arbitrariness of the function representation as local
gauge freedom:
\begin{definition}
  \begin{enumerate}
    \item Let $(M,H,E)$ be an operator manifold.
        An isomorphism
\Equ{ZZ}
   V \::\: H \rightarrow \bigoplus_{m \in \sN \cup \{0, \infty\}} (L^2(D_m,\mu))^m
\EndEqu
        is called a {\bf{gauge}}.
    \item For two gauges $V$,$\tilde{V}$, the mapping
\[        U = V \tilde{V}^{-1} : \bigoplus_{m \in \sN \cup \{0, \infty\}}
	  (L^2(D_m,\mu))^m \rightarrow \bigoplus_{m \in \sN \cup \{0, \infty\}}
	  (L^2(D_m,\mu))^m          \]
        is called a {\bf{gauge transformation}}.
        It can, according to \Ref{loc_unit1},\Ref{repr}, be represented as local
	$U(m_x)$-transformation of the wave functions.
  \end{enumerate}
\end{definition}
The occurrence of local gauge freedom can, on a non-rigorous level, be
understood from the fact that only $|\Psi(x)|^2 = \sum_{\alpha=1}^{m_x}
|\Psi^\alpha(x)|^2$ is a physical observable. The local gauge group $U(m_x)$
is the isometry group of the spin scalar product.

If taken seriously, our concept has important physical consequences:
The local gauge principle is no longer a fundamental physical principle,
but follows from the fact that space is a quantum mechanical observable.
In contrast to usual gauge theories, the gauge group cannot be chosen arbitrarily.
For a given configuration of the spinors, it is fixed to be the group $U(m_x)$, which
acts directly on the spinorial index of the wave functions. This is a strong
restriction for the formulation of physical models.

We point out that the transformation functions $W^\alpha_\beta$ in \Ref{repr}
are in general not smooth, they are only measurable.
From the physical point of view, it seems reasonable to restrict to smooth
gauge transformations. Then the structure of an operator manifold
reduces to a principle bundle over $M$ with fibre $\:\C^m$ and local gauge
group $U(m)$. The wave functions are sections of the bundle.
In this way, we obtain the mathematical framework of classical gauge theories.
In view of the fact that the restriction to smooth gauge transformations is only
a technical convenience, however, the question arises if the topology of the
fibre bundles has physical significance.

We remark that our constructions have a direct generalization to
relativistic quantum mechanics \cite{F1}. The adaptation to many-fermion systems
finally leads to the ``Principle of the Fermionic Projector'' as introduced
in \cite{F2}.

\section{\bf Acknowledgments.}
I thank Claus Gerhardt for his interest and support.

\addcontentsline{toc}{section}{References}

\nonumsection{References}

\end{document}